\documentclass[12pt]{amsart}
\usepackage{amsfonts}
\usepackage{amssymb}
\theoremstyle{definition}

\theoremstyle{remark}

\begin{document}

\centerline{\large\bf THE AXIOM OF COHOLOMORPHIC $(2n+1)$-SPHERES}
\centerline{\large\bf IN THE ALMOST HERMITIAN GEOMETRY
\footnote{\it SERDICA Bulgaricae mathematicae publicationes. Vol. 8, 1982, p. 391-394}}

\vspace{0.2in}
\centerline{\large OGNIAN  KASSABOV}

\vspace{0.5in}
{\sl In his book on Riemannian geometry [1] E. C\,a\,r\,t\,a\,n proved a characterization
of a real-space-form, using the axiom of planes. There are many results in this direction
also for a Kaehler manifold. B.-Y. C\,h\,e\,n and K. O\,g\,i\,u\,e [4] have proved that a
Kaehler manifold, which satisfies the axiom of coholomorphic 3-spheres is flat. In this 
paper we prove a generalization of this theorem for an almost Hermitian manifold.}

\vspace{0.3in}
{\bf 1. Introduction.} Let $N$ be an $n$-dimensional submanifold of an $2m$-dimensional
almost Hermitian manifold $M$ with Riemannian metric $g$ and almost complex structure $J$.
Let $\tilde\nabla$ and $\nabla$ be the covariant differentiations on $M$ and $N$, 
respectively. It is well known, that the equation 
$
	\alpha(X,Y) = \widetilde\nabla_XY-\nabla_XY  ,
$
where $X,\ Y \in {\mathfrak X}N$, defines a normal-bundle-valued symmetric tensor field, called
the second fundamental form of the immertion. The submanifold $N$ is said to be totally
umbilical, if $\alpha(X,Y)=g(X,Y)H$ for all $X,\ Y \in {\mathfrak X}N$, where $H=(1/n){\rm trace}\, \alpha$
\ is the mean curvature vector of $N$ in $M$. In particular, if $\alpha$ vanishes identically,
$N$ is called a totally geodesic submanifold of $M$.

For $X \in \mathfrak XN,\ \xi \in \mathfrak XN^{\perp}$ we write 
$ \widetilde\nabla_X\xi =- A_{\xi}X+D_X{\xi}$, where $-A_{\xi}X$ (respectively, $D_X\xi$)
denotes the tangential (resp. the normal) component of $\widetilde\nabla_X\xi$. A normal
vector field $\xi$ is said to be parallel, if $D_X\xi=0$ for each $X\in \mathfrak XN$.

By an $n$-plane we mean an $n$-dimensional linear subspace of a tangent space. 
A $2n$-plane (respectively an $n$-plane) \ $\alpha$ \ where $1\le n\le m $ \ is said to be 
holomorphic (respectively, antiholomorphic) if  \ $J\alpha = \alpha$ \ 
(respectively \ $J\alpha \perp \alpha)$. A $(2n+1)$-plane $\alpha$ is called coholomorphic 
if it contains a holomorphic $2n$-plane. 

An almost Hermitian manifold $M$ is said to satisfy the axiom of holomorphic $2n$-planes 
(respectively $2n$-spheres) if for each point $p \in M$ and for any $2n$-dimensional
holomorphic plane \ $\pi$ \ in \ $T_pM$ \ there exists a totally geodesic
submanifold \ $N$ \ (respectively a totally umbilical submanifold \ $N$ \ with nonzero
parallel mean curvature vector) containing \ $p$, such that \ $T_pN=\pi$, where $n$ is
a fixed number, $1\le n<m$.

An almost Hermitian manifold $M$ is said to satisfy the axiom of antiholomorphic $n$-planes 
(respectively $n$-spheres) if for each point $p \in M$ and for any $n$-dimensional
antiholomorphic plane \ $\pi$ \ in \ $T_pM$ \ there exists a totally geodesic
submanifold \ $N$ \ (respectively a totally umbilical submanifold \ $N$ \ with nonzero
parallel mean curvature vector) containing \ $p$, such that \ $T_pN=\pi$, where $n$ is
a fixed number, $1< n\le m$.

An almost Hermitian manifold is called an \ $RK$-manifold, if  \ 
$
	R(X,Y,Z,U)=R(JX,JY, \\ JZ,JU)
$ 
for all \ $X,\, Y,\, Z,\, U \in T_pM$, \ $p\in M$. 

We have proved in [5]:

\vspace{0.05in}
T\,h\,e\,o\,r\,e\,m A. {\it Let \ $M$ \ be a \ $2m$-dimensional almost Hermitian manifold,
$m \ge 2$. If \ $M$ \ satisfies the axiom of holomorphic \ $2n$-planes or the axiom of
holomorphic \ $2n$-spheres for some \ $n$, $1\le n <m$, then \ $M$ \ is an \ 
$RK$-manifold with pointwise constant holomorphic sectional curvature.}

\vspace{0.05in}
T\,h\,e\,o\,r\,e\,m B. {\it Let  $M$  be a  $2m$-dimensional almost Hermitian manifold,
$m \ge 2$. If  $M$  satisfies the axiom of antiholomorphic \ $n$-planes or the axiom of
antiholomorphic  $n$-spheres for some  $n$, $1< n \le m$, then  $M$  is an  
$RK$-manifold with pointwise constant holomorphic sectional curvature and with pointwise
constant antiholomorphic sectional curvature. Consequently, if $m\ge 3$, then \ $M$ \ is one 
of the following:

1) a real-space-form, or

2) a complex-space-form.}

\vspace{0.05in}
These theorems generalize some results in [3, 6, 9].
It is not difficult to see that if $n>1$ then the holomorphic analogue of Theorem B holds.

Following B.-Y. C\,h\,e\,n and K. O\,g\,i\,u\,e [4], L. V\,a\,n\,h\,e\,c\,k\,e formulates
the following axiom of coholomorphic $(2n+1)$-spheres [8]:

\vspace{0.05in}
{\it For each point $p \in M$ and for each coholomorphic $(2n+1)$-plane $\pi$ in $T_pM$,
there exists a $(2n+1)$-dimensional totally umbilical submanifold $N$ of $M$ containing
$p$, such that $T_pN=\pi$, where $n$ is a fixed integer, $1\le n <m$.}

\vspace{0.05in}
We shall prove the following theorem.

\vspace{0.05in}
T\,h\,e\,o\,r\,e\,m. {\it Let  $M$  be a  $2m$-dimensional almost Hermitian manifold,
$m \ge 2$. If  $M$  satisfies the axiom of coholomorphic  $(2n+1)$-spheres for 
some  $n$, then  $M$  is conformal flat.}

\vspace{0.05in}
Hence, using [7] we have

\vspace{0.05in}
C\,o\,r\,o\,l\,l\,a\,r\,y 1. {\it Let $M$ be a $2m$-dimensional connected Kaehler
manifold, $m\ge 2$. If  $M$  satisfies the axiom of coholomorphic  $(2n+1)$-spheres for 
some  $n$, then either $M$ is flat or $M$ is locally a product of two 2-dimensional
Kaehler manifolds with constant curvature $K$ and $-K$, respectively, $K>0$.}

\vspace{0.05in}
The case $m\ge 3$ in corollary 1 is treated in [4].

An almost Hermitian manifold $M$ which satisfies $(\widetilde\nabla_XJ)X=0$ for all
$X\in \mathfrak XM$ is said to be an $NK$-manifold. Using the classification in [7]
we have also

\vspace{0.05in}
C\,o\,r\,o\,l\,l\,a\,r\,y 2. {\it Let $M$ be a $2m$-dimensional $NK$-
manifold, $m\ge 2$. If  $M$  satisfies the axiom of coholomorphic  $(2n+1)$-spheres for 
some  $n$, then  $M$ is one of the following:

1) a flat Kaehler manifold,

2) locally a product $M_1\times M_2$, where $M_1$ (respectively $M_2$) is a 
2-dimensional Kaehler manifold with constant curvature $K$ (respectively $-K$),

3) a 6-dimensional manifold of constant curvature $K>0$,

4) locally a product  $M_3\times M_2$, where $M_3$  is a 
6-dimensional $NK$-manifold of constant curvature $K>0$.}

\vspace{0.05in}
An almost Hermitian manifold $M$ is said to be of pointwise constant type $\alpha$,
provided that for each point $p\in M$ and for each $X\in T_pM$ we have 
$\alpha(p)g(X,X)=\lambda(X,Y)=\lambda(X,Z)$ with 
$\lambda(X,Y)=R(X,Y,Y,X)-R(X,Y,JY,JX)$ whenever the planes defined by $X,Y$ and
$X,Z$ are antiholomorphic and $g(Y,Y)=g(Z,Z)=1$. If for $X,Y \in {\mathfrak X}(M)$
with $g(JX,Y)=g(X,Y)=0$, $\lambda(X,Y)$ is a constant whenever $g(X,X)=g(Y,Y)=1$,
then $M$ is said to have global constant type.

\vspace{0.05in}
C\,o\,r\,o\,l\,l\,a\,r\,y 3. {\it Let $M$ be an almost Hermitian manifold with pointwise
constant type $\alpha$. If  $M$  satisfies the axiom of coholomorphic  $(2n+1)$-spheres for 
some  $n$ and if ${\rm dim} M \ge 6$, then $M$ is a space of constant curvature 
$\alpha$ and $M$ has global constant type.}

\vspace{0.05in}
Corollary 3 is proved in [8] for an $RK$-manifold.

\vspace{0.2 in}
{\bf 2. Preliminaries.} Let $M$ be an $2m$-dimensional almost Hermitian manifold with
Riemannian metric $g$, almost complex structure $J$ and covariant
differentiation $\nabla$. The curvature tensor $R$, associated with $\nabla$ 
has the following properties:

1) $ R(X,Y)Z=-R(Y,X)Z$

2) $ R(X,Y)Z+R(Y,Z)X+R(Z,X)Y=0$ 

3) $ R(X,Y,Z,U)=-R(X,Y,U,Z)$

\noindent for all $X,\ Y,\ Z,\ U \in T_pM$, $p\in M$, where $ R(X,Y,Z,U)=g(R(X,Y)Z,U)$. 
The Weil conformal curvature tensor $C$ is defined by
$$
	\begin{array}{r}
		C(X,Y,Z,U)=R(X,Y,Z,U)-(1/(2m-2))\{ g(X,U)S(Y,Z)  \\
		-g(X,Z)S(Y,U)+g(Y,Z)S(X,U)-g(Y,U)S(X,Z) \}   \\
		+(S(p)/((2m-1)(2m-2))) \{ g(X,U)g(Y,Z)-g(X,Z)g(Y,U) \},
	\end{array}
$$
where $S$ and $S(p)$ are the Ricci tensor and the scalar curvature of $M$, respectively.

Now, let $N$ be a submanifold of $M$, as in section 1. The normal component of
$R(X,Y)Z$, where $X,Y,Z \in {\mathfrak X}N$ is given by
$$
	(R(X,Y)Z)^{\perp} = (\overline\nabla_X\alpha)(Y,Z)-(\overline\nabla_Y\alpha)(X,Z) \ ,  \leqno (2.1)
$$
where $(\overline\nabla_X\alpha)(Y,Z)=D_X\alpha(Y,Z)-\alpha(\nabla_XY,Z)-\alpha(Y,\nabla_XZ) $
and if $N$ is totally umbilical submanifold of $M$, (2.1) reduces to
$$
	(R(X,Y)Z)^{\perp} = g(Y,Z)D_XH-g(X,Z)D_YH \ .  \leqno (2.2)
$$

\vspace{.1in}
{\bf 3. Proof of the Theorem}. Let $X,\,Y$ be arbitrary unit vectors in $T_pM$,
$p\in M$, such that $X$ is perpendicular to $Y,\,JY$. Applying the axiom of coholomorphic 
$(2n+1)$-spheres for a coholomorphic plane, which contains $X,\,JX,\,JY$ and is
perpenducular to $Y$ and using (2.2) we obtain
$$
	R(X,JX,JY,Y)=0 \ ,          \leqno(3.1)
$$
$$
	R(JY,JX,X,Y)=0 \ ,         \leqno (3.2)
$$
$$
	R(X,JX,JX,Y)=g(D_XH,Y) \ ,
$$
$$
	R(X,JY,JY,Y)=g(D_XH,Y) \ .
$$
Hence
$$
	R(X,JX,JX,Y)=R(X,JY,JY,Y) \ .   \leqno (3.3)
$$

From (3.2) we have
$ R(Y+JY,JX,X,Y-JY)=0$ and consequently
$$
	R(X,Y,Y,JX)=R(X,JY,JY,JX) \ .   \leqno (3.4)
$$

If $m>2$, we take a unit vector $Z$, perpenducular to $X,\,JX,\,Y,\,JY$.
Using again the axiom of coholomorphic $(2n+1)$-spheres and (2.2) we find
$$
	R(X,JX,Y,Z)=R(X,Y,JY,Z)=0 \ ,       \leqno (3.5)
$$
$$
	R(X,JX,JX,Z)=R(X,Y,Y,Z) \ ,       \leqno (3.6)
$$
$$
	R(X,Y,Y,JX)=R(X,Z,Z,JX) \ .       \leqno (3.7)
$$

If $m\ge 4$, let $U$ be a unit vector in $T_pM$, perpenducular to
$X,\,JX,\,Y,\,JY,\,Z,\,JZ$. From (3.6) we have
$ 2R(X,JX,JX,U)=R(X,Y+Z,Y+Z,U)=0$, which gives
$R(X,Y,Z,U)=-R(X,Z,Y,U)$.

Hence, by the properties of the curvature tensor $R$ we obtain
$$
	R(X,Y,Z,U)=0  \ .           \leqno (3.8)   
$$

Making use of (3.1)-(3.8) it is not difficult to prove that $R(X,Y,Z,U)=0$
for an arbitrary orthogonal quadriple $X,\,Y,\,Z,\,U \in T_pM$. According
to a well known theorem of Schouten [2] the Weil conformal curvature 
tensor of $M$ vanishes.

R\,e\,m\,a\,r\,k. If a Riemannian manifold $M$ of dimension $n>3$ is 
conformal flat, then there exists a totally umbilical submanifold $N$ of
dimension $n<m$ through every point of $M$ and in every $n$-dimensional
direction of that point (see [2]). Consequently, if $M$ is a conformal flat
$2m$-dimensional almost Hermitian manifold, $m\ge 2$, then $M$ satisfies
the axiom of coholomorphic $(2n+1)$-spheres for every $n,\,1\le n<m$.

\vspace{0.4in}
\centerline{\large R E F E R E N C E S}

\vspace{0.1in}
\noindent
1. E. C\,a\,r\,t\,a\,n. Le\c cons sur la g\'eometrie des espaces de Riemann. Paris, 1946.

\noindent
2. J. S\,c\,h\,o\,u\,t\,e\,n. Ricci-calculus. Berlin, 1954.

\noindent
3. B. -Y. C\,h\,e\,n, K. O\,g\,i\,u\,e. Some characterizations of complex space forms. {\it Duke

 Math. J.}, {\bf 40}, 1973, 797-799.

\noindent
4. B. -Y. C\,h\,e\,n, K. O\,g\,i\,u\,e. Two theorems on Kaehler manifolds. {\it Michigan Math. J.,}

{\bf 21}, 1974, 225-229.
 
\noindent
5. O. K\,a\,s\,s\,a\,b\,o\,v. On the axiom of planes and the axiom of spheres in the almost Hermi-

tian geometry.
{\it Serdica}, {\bf 8}, 1982, 109-114.

\noindent
6. K. N\,o\,m\,i\,z\,u. Conditions for constancy of the holomorphic sectional curvature.

{\it J. Diff. Geom.}, {\bf 8}, 1973, 335-339.

\noindent
7. S. T\,a\,n\,n\,o. 4-dimensional conformally flat Kaehler manifolds. {\it Tohoku Math. J.}, {\bf 24}, 

1972, 501-504. 

\noindent
8. L. V\,a\,n\,h\,e\,c\,k\,e. The axiom of coholomorphic (2p+1)-spheres for some almost Hermit-

ian manifolds. 
{\it Tensor (N.S.)}, {\bf 30}, 1976, 275-281.

\noindent
9. K. Y\,a\,n\,o, I. M\,o\,g\,i. On real representation of Kaehler manifolds. {\it Ann. Math.}, {\bf 61,}

 1955, 170-189.

\vspace {0.3cm}
\noindent
{\it Center for mathematics and mechanics \ \ \ \ \ \ \ \ \ \ \ \ \ \ \ \ \ \ \ \ \ \ \ \ \ \ \ \ \ \ \ \ \ \
Received 18.12.1980

\noindent
1090 Sofia   \ \ \ \ \ \ \ \ \ \ \ \ \ \ \ \ \  P. O. Box 373}

\end{document}